\documentclass{article}
\usepackage{amssymb}
\usepackage{latexsym}
\newtheorem{thm}{Theorem}
\newtheorem{tom}{Theorem}
\newtheorem{lem}{Lemma}
\newtheorem{cor}{Corollary}

\begin{document}
\title{Minimal Lagrangian tori in Kahler-Einstein manifolds}
\author{Edward Goldstein}
\maketitle

\renewcommand{\abstractname}{Abstract}
\begin{abstract}
In this paper we use structure preserving torus actions on Kahler-Einstein
manifolds to construct minimal Lagrangian submanifolds. Our main result is: Let
$N^{2n}$ be a Kahler-Einstein manifold with positive scalar curvature with an
effective $T^n$-action. Then precisely one regular orbit $L$ of the $T$-action
is a minimal Lagrangian submanifold of $N$. Moreover there is an $(n-1)$-torus
$T^{n-1} \subset T^n$ and a sequence of non-flat immersed minimal Lagrangian
tori $L_k$, invariant under $T^{n-1}$ s.t. $L_k$ locally converge to $L$ (in particular the supremum of the sectional curvatures of
$L_k$ and the distance between $L$ and $L_k$ go to $0$ as $k \mapsto \infty$).
\end{abstract}
\section{Introduction}
In this paper we will use torus actions on Kahler-Einstein manifolds with
positive scalar curvature to construct minimal Lagrangian tori.

Let $N^{2n}$ be a Kahler-Einstein (K-E) manifold with positive scalar
curvature and suppose we have a structure preserving $T^k$-action on $N$.
We will look for $T$-invariant minimal Lagrangian submanifolds of $N$.
If $k=n$ then we have shown in \cite{Gold} that there is precisely one regular
orbit of the $T$-action, which is a minimal Lagrangian submanifold of $N$. In
this paper we study the case when $k=n-1$ (a complexity one action).
The main tool in our investigation will be a correspondence between minimal
Lagrangian submanifolds of $N$ and certain Special Lagrangian submanifolds of
$K(N)$- the total space of the canonical bundle of $N$. 

The manifold $K(N)$ has a natural holomorphic volume form $\varphi$.
Also since $N$ is K-E with positive scalar curvature, we have a (Calabi)
metric $\omega_u$ on $K(N)$, which is a Ricci-flat Kahler metric
(see Section 2.1). The form $\varphi$ is covariantly constant with respect
to $\omega_u$, and we have Special Lagrangian (SLag) submanifolds
$L' \subset K(N)$, defined by the conditions $\omega_u|_{L'}=0$ and
$Im\varphi |_{L'}=0$ (see \cite{HL} and \cite{Gold}). There is a radial vector
field $Y$ on $K(N)$, whose flow is scaling of $K(N)$ by real numbers
(see Section 2.1). Our main tool in studying minimal Lagrangian submanifolds
on $N$ will be a correspondence between minimal Lagrangian submanifolds on $N$
and SLag submanifolds on $K(N)$, invariant under the flow of $Y$
(see Lemmas 1 and 2 in Section 2.1).

In Section 2.2 we study SLag submanifolds on $K(N)$ using a torus action
on $N$. Suppose we have a $T^k$-action on $N$. This action of course induces a
$T^k$-action on $K(N)$. We will see that there are canonical moment maps $\mu$
on $N$ and $\mu'$ on $K(N)$. Let $Z \subset N$ be the zero set of $\mu$ and
$\pi:K(N) \mapsto N$ be the projection. Then the zero set of $\mu'$ is
$Z'=\pi^{-1}(Z)$. Suppose that $T$ acts freely on $Z''=Z'-Z$. Then we have
a symplectic reduction $Q=Z''/T$. We will see
that $Q$ has a natural holomorphic volume form $\varphi'$ and a metric
$\omega'$ and SLag submanifolds of $(Q,\varphi', \omega')$ lift
to $T$-invariant SLag submanifolds of $K(N)$. Also the vector field $Y$ is
tangent to $Z''$ and projects to a vector field $Y'$ on $Q$. Thus
we reduced the problem of finding minimal Lagrangian submanifolds of $N$ to a
problem of finding SLag submanifolds of $Q$, invariant under the flow of
$Y'$.

In Section 2.3 we assume that $k=n-1$. Let $X \subset Z''$ be the
set of elements of $Z''$ of unit length in $K(N)$ and $S=X/T \subset Q$. We
will see that there is a non-vanishing vector field $W$ on $S$
s.t. there is a correspondence between $Y'$-invariant SLag submanifolds of
$Q$ and trajectories of the $W$-flow on $S$.

Next we would like to develop a criterion to see that $T^{n-1}$ acts freely
on $Z''$. We also would like to understand periodic orbits of the vector
field $W$ on $S$ (to construct immersed minimal Lagrangian tori on $N$).
We can do it if we assume that $N$ is a toric K-E manifold (see Section 3.1).
In this case we can prove the following Theorem: 
\begin{tom}
Let $N^{2n}$ be a K-E manifold with positive scalar curvature with an
effective $T^n$-action. Then precisely one regular orbit $L$ of the $T$-action
is a minimal Lagrangian submanifold of $N$.
Moreover there is an $(n-1)$-torus $T^{n-1} \subset T^n$ and a sequence of
non-flat immersed $T^{n-1}$-invariant minimal Lagrangian tori
$L_k \subset N$ s.t. $L_k$ locally converge to $L$ (in particular the supremum of sectional curvatures of $L_k$ and the distance
between $L$ and $L_k$ go to $0$ as $k \mapsto \infty$). 
\end{tom}
Here by local convergence we mean the following: The distance between $L_k$ and $L$ goes to $0$ as $k \mapsto \infty$. Also for any point $l \in L$ we can choose a neighbourhood $U$ of $l$ in $N$ s.t. $L_k \bigcap U$ is a finite union $L_k^j$ of submanifolds of the form $L_k^j=exp(v_k^j)(L \bigcap U)$, where $v_k^j$ is a normal vector field to $L$ on $L \bigcap U$. Moreover any subsequence $v_k^j$ converges to $0$ in a $C^{\infty}$ topology as $k \mapsto \infty$.

This result is new even for $N= \mathbb{C}P^n$ for $n \geq 3$. For $n=2$ examples of non-flat $S^1$-invariant immersed minimal Lagrangian tori in $\mathbb{C}P^2$ were constructed in \cite{Cas} and \cite{Has} using harmonic maps. 

{\bf Acknowledgments:} This paper is a part of author's work towards his Ph.D. at MIT. He wants to thank his advisor, Tom Mrowka, for continuing support. 

In this paper we use a number of results from our previous paper \cite{Gold}, including proofs for the completeness of exposition.

Research is partially supported by an NSERC PGS B award. 
\section{Minimal Lagrangian submanifolds in complexity one K-E manifolds}
\subsection{A correspondence between minimal and Special Lagrangian
submanifolds}
Let $N^{2n}$ be a K-E manifold with positive scalar curvature. We begin by reviewing the geometry of $K(N)$ and the correspondence between minimal Lagrangian submanifolds of $N$ and certain Special Lagrangian submanifolds of $K(N)$.

Let $K(N)$ be the total space of the canonical bundle of $N$ and $\pi : K(N) \mapsto N$ be the projection. There is a canonical $(n,0)$-form $\rho$ on $K(N)$ defined by $\rho(a)(v_1,\ldots,v_n)=a(\pi_{\ast}(v_1),\ldots,\pi_{\ast}(v_n))$, $a \in K(N)$.
The form $\varphi = d\rho$ is a holomorphic volume form on $K(N)$. If $z_1, \ldots,z_n$ are local coordinates on $N$ then
 $(z_1,\ldots,z_n,y=dz_1\wedge \ldots \wedge dz_n)$ are coordinates on $K(N)$ and $\rho=ydz_1 \wedge \ldots \wedge dz_n$, $\varphi= dy \wedge dz_1\wedge \ldots \wedge dz_n$.

There is a canonical radial vector field $Y$ on $K(N)$, given at a point $m \in K(N)$ by the vector $m$ (viewed as a tangent vector to the linear fiber over $\pi(m)$). We have $i_Y\rho=0$. Also the Lie derivative ${\cal L}_Y \rho= \rho$. So $\rho= i_Y d\rho=i_Y \varphi$. So ${\cal L}_Y \varphi= d(i_Y \varphi)=d\rho=\varphi$. 

If $N$ is a Kahler-Einstein manifold with positive scalar curvature then $K(N)$ has a
Ricci-flat Kahler metric on it (see \cite{Sol}, p.108). The metric is constructed as follows : The connection on $K(N)$ induces a horizontal distribution for the projection $\pi$, with a corresponding splitting of the tangent bundle of $K(N)$ into horizontal and vertical distributions. We can identify the horizontal space at each point $m \in K(N)$ with the tangent space to $N$ at $\pi(m)$. Let $r^2 : K(N) \mapsto 
\mathbb{R}_+$ be the square of the length of an element in $K(N)$ and $u: \mathbb{R}_+ \mapsto \mathbb{R}_+$ be a positive function with a positive first derivative. We define the metric $\omega_u$ on $K(N)$ as follows: We put the horizontal and the vertical distributions to be orthogonal. On the horizontal distribution we define the metric to be $u(r^2)\pi^{\ast}(\omega)$ and on the vertical distribution we define it to be $t^{-1}u'(r^2)\omega^{\star}$. Here $\omega$ is the Kahler-Einstein metric on $N$, $t$ is its scalar curvature and $\omega^{\star}$ is the induced metric on the linear fibers of $\pi $. The Kahler-Einstein condition on $N$ ensures that the corresponding 2-form $\omega_u$ defining this metric on $K(N)$ is closed, i.e. the metric is Kahler. If we take $u(r^2)=(tr^2+l)^{\frac{1}{n+1}}$ for some positive constant $l$  (see \cite{Sol}, p.109), then $\omega_u$ is complete and Ricci-flat (the Calabi metric). From now on we study $K(N)$ endowed with this metric $\omega_u$.

We begin with the following observation : Let $L$ be an oriented Lagrangian submanifold of $N$. For any point $l \in L$ there is
a unique element $\kappa_l$ in the fiber of $K(N)$ over $l$ which restricts to the volume form on $L$. Various $\kappa_l$ give rise to a section $\kappa$ of $K(N)$ over $L$. Consider a submanifold $L^K \subset K(N)$ given by \[ L^K=((m|
m=a\kappa_l ~ for ~ l \in L ~ , ~ a \in \mathbb{R}) \] 
We have the following:
\begin{lem}
$L$ is a minimal Lagrangian submanifold of $N$ iff $L^K$ is a Special Lagrangian submanifold of $K(N)$
\end{lem} 
Here by a minimal submanifold we mean a submanifold, which is critical for the volume functional (i.e. the trace of the second fundamental form vanishes).

{\bf Proof :} First we note that $L^K$ is Special, i.e.
$Im\varphi|_{L^K}=0$. Indeed one easily verifies that $Im\rho|_{L^K}=0$, hence $Im\varphi|_{L^K}=0$.

We now prove that $L^K$ is Lagrangian with respect to $\omega_u$ iff $L$ is minimal. Let $m$ be a point on $L^K-L$, $l = \pi(m)$ and $m = a\kappa_l$. The tangent space of $L^K$ at $m$ is spanned by $\kappa_l$ (viewed as a vertical vector in $T_mK(N)$) and vectors $(e+a
\nabla_e\kappa)$. Here $e$ is any tangent vector to $L$ at $l$ (viewed as an element of the horizontal
distribution of $T_mK(N)$) and $a\nabla_e \kappa$ lives in the vertical distribution of $T_mK(N)$. To compute $\nabla_e \kappa$ take an orthonormal frame $(v_j)$ of $T_lL$ and extend it to an orthonormal frame of $L$ in a neighbourhood $U$ of $l$ in $L$ s.t. $\nabla^L v_i=0$ at $l$ (here $\nabla^L$ is the Levi-Civita connection of $L$). We get that \[ \nabla_e \kappa = \kappa \cdot \nabla_e \kappa(v_1,\ldots,v_n)= \kappa (e(\kappa(v_1,\ldots,v_n))- \Sigma\kappa(v_1,\ldots,\nabla_e v_j,\ldots,v_n)) \]
Now $e(\kappa(v_1,\ldots,v_n))=0$. Also clearly \[\kappa(v_1,\ldots,\nabla_e v_j,\ldots,v_n)=i<\nabla_e v_j,Jv_j>=i<\nabla_{v_j}e,Jv_j>=i<-e,J(\nabla_{v_j}v_j)> \]   
Here $J$ is the complex structure on $N$. Thus we get that 
\[a \nabla_e \kappa = -ia(Jh \cdot  e)\kappa_l \]
Here $h= \Sigma \nabla_{v_j}v_j$ is the trace of the second fundamental form of $L$. From this one easily deduces that $L^K$ is Lagrangian iff $h=0$, i.e. $L$ is minimal.    Q.E.D.

The manifold $L^K$ is invariant under the flow of the vector field $Y$ on $K(N)$ (which is just scaling of $K(N)$ by real numbers). Vice versa we have the following:
\begin{lem}
Let $L'$ be a Special Lagrangian submanifold of $K(N)-N$, invariant under the flow of $Y$. Then $L=\pi(L')$ is an (immersed) minimal Lagrangian submanifold of $N$.
\end{lem}
{\bf Proof:} Let $m \in L'-N$. Since $L'$ is Lagrangian and $Y$ is in the tangent space $T_m L'$ then the tangent space to $L'$ at $m$ clearly decomposes as
\[T_m L'= span(Y) \oplus T'\] where $T'$ is in the horizontal distribution at $m$. The space $\pi_{\ast}(T')$ can be viewed as a tangent space to $L$ at $l= \pi(m)$. Clearly this tangent space $T_l L$ is Lagrangian, i.e. $L$ is Lagrangian. Also $L'$ was Special and we have seen in the beginning of this section that $i_Y \varphi= \rho$. Thus $m$ (viewed as an $(n,0)$-form on $N$ at $l$) restricts to a real $n$-form on $T_l L$, i.e. $m \in L^K$. Hence locally 
$L'$ coincides with $L^K$. From Lemma 1 we deduce that $L$ is minimal. Q.E.D.

\subsection{SLag submanifolds on $K(N)$ via symplectic reduction}
 
In the previous section we showed how to find minimal Lagrangian submanifolds
of $N$ from certain SLag submanifolds of $K(N)$. In this section we will see
that if we have a torus action on $N$ then we can find $T$-invariant SLag
submanifolds of $K(N)$ from SLag submanifolds of a certain symplectic
reduction of $K(N)$.

Let $T^k$ act on $N$. Then this action induces a $T^k$-action on $K(N)$.
Let $\mathcal{T}$ be the Lie algebra of $T$, $ v \in \mathcal{T}$, $X_v$ be
the flow vector field on $N$ and $X_v'$ the flow vector field on $K(N)$.
So $\pi_{\ast}(X_v')=X_v$. Let $l \in N$ and $m \in K_l=\pi^{-1}(l)$.
Let $R(m)$ be the vertical part of $X_v'$ at $m$. Since $R(m)$ is vertical, it
can be viewed as an element of $K_l$. The correspondence $m \mapsto R(m)$ is a
linear correspondence on $K_l$. Hence there is a complex number $\sigma_l(v)$
s.t. $R(m)=\sigma_l(v)m$. At a regular point $l$ of the $T$-action $\sigma_l(v)$ can also be found in a following way : Take any unit length element $\xi \in K_l$. Extend $\xi$ along the orbit of $X_v$ to be invariant under the flow of $X_v$. Then one easily computes that $\sigma_l(v)= \nabla_{X_v}\xi \cdot \xi$. Since the flow of $X_v$ is given by holomorphic isometries, $\xi$ has unit length. Hence  $\sigma_l(v)$ is purely imaginary. Also $\sigma_l(v)$ is linear in $v$ (because $R(m)$ is given by the vertical part of the differential of the $T$-action at $m$, and this differential is a linear map from $\mathcal{T}$ to $T_mK(N)$). Hence $i \sigma$ can be viewed as a map from $N$ to the dual Lie algebra $\mathcal{T}^{\ast}$. This map is $T$-invariant.
 
Let $t > 0$ be the scalar curvature of $N$.
\begin{lem}
The map $\mu = -it^{-1} \sigma$ is a moment map for the action.
\end{lem}
{\bf Proof:} Let $v \in \mathcal{T}$. We need to show that $d(-it^{-1}\sigma(v))= i_{X_v}\omega$. We will do it at a regular point $l$ of the action. Choose any unit length element $\xi$ of $K(N)$ over $l$. We can extend $\xi$ to be a local unit length section, invariant under the $X_v$-flow. $\xi$ defines a connection 1-form $\psi$, $\psi(u)= \nabla_u \xi \cdot \xi$. $\psi$ is invariant under the $X_v$-flow and the K-E condition says that $i d\psi = t \omega$. 
So \[ 0 = \mathcal{L}_{X_v}\psi= d(i_{X_v}\psi) + i_{X_v}d\psi= d\sigma(v) - it (i_{X_v}\omega) \]   
So $\mu$ is a moment map.   Q.E.D. 

{\bf Remark:} By the construction of $\mu$ we get that $\mu(v)=0$ for some $v \in {\cal T}$ at a point $l \in N$ iff the vector field $X_v'$ is horizontal at $\pi^{-1}(l)$.
\begin{lem} 
The map $\mu'= u\pi^{-1}(\mu)$ is a moment map for the $T$-action on
$K(N)$.
\end{lem}
{\bf Proof:} Let $v \in \cal{T}$. We need to prove that $d\mu'(v)=i_{X_v'}\omega_u$.
 
We will study $\omega_u$ in more detail (see \cite{Sol}). Let $m \in N$ be a regular point for the $T^n$-action and $\xi$ a unit length element of $K(N)$ over $m$. We can extend $\xi$ to be a local unit length section of $K(N)$, invariant under the flow of $X_v$. $\xi$ gives rise to a connection 1-form $\psi$ for the connection on $K(N)$  and the Einstein condition tells that $id\psi= t \omega$. The section $\xi$ defines a complex coordinate $a$ on $K(N)$, which is invariant under the $X_v'$-flow. Also the form $b=da + a\pi^{\ast}\psi$ vanishes on the horizontal distribution (see \cite{Sol}, p. 108). We have $r^2=a\overline{a}$ and $u= u(r^2)$. Also the Kahler form $\omega_u$ on $K(N)$ is given by \[\omega_u=u\pi^{\ast}\omega -it^{-1}u'b\wedge \overline{b} \]
One directly verifies that $\omega_u = d\eta$ for $\eta= it^{-1}u \pi^{\ast} \psi - it^{-1}\frac{ud\overline{a}}{\overline{a}}$. By our construction the flow of $X_v'$ leaves $\eta$ invariant. So \[0 = {\cal L}_{X_v'}\eta= i_{X_v'}d\eta + d(i_{X_v'}\eta)= i_{X_v'}\omega_u + d(it^{-1}u\psi(X_v))= i_{X_v'}\omega- d(\mu'(v))\] Here we used the fact that $d\overline{a}(X_v')= 0$ and $\psi(X_v)= \sigma(v)$. So $\mu'$ is a moment map and we are done.   Q.E.D.

Let now $L'$ be a (connected) SLag submanifold of $K(N)$,
invariant under the $T$-action and under the $Y$-flow.
Since $L'$ is Lagrangian and $T$-invariant, the moment map $\mu'$ is constant
on $L'$. But $\mu'= u \pi^{-1}(\mu)$ and $Y(\mu')= 2r^2 u' \pi^{-1}(\mu)$. So
we have $\pi^{-1}(\mu)=0$ on $L'$. Let $Z$ be the zero set of $\mu$. Then
$L' \subset Z' = \pi^{-1}(Z)=\mu'^{-1}(0)$. 

Let $Z''=Z'-Z$. From now on we assume that $T$ acts freely on $Z''$ (we
will demonstrate examples where this holds in Section 3). 
We have a symplectic reduction $N_{red}=Z/T$ and a (smooth)
symplectic reduction $Q=Z''/T$, endowed with a Kahler metric $\omega'$.

Let $v_1, \ldots, v_k$ be a basis for $\mathcal{T}$ and $X_1', \ldots, X_k'$
be the corresponding flow vector fields on $K(N)$.
Let $\varphi^{\ast}=i_{X_1'} \ldots i_{X_k'}\varphi$ be an $(n-k,0)$-form
on $K(N)$, obtained by contracting $\varphi$ by $X_1',\ldots,X_k'$.
Let $\rho^{\ast}=i_{X_1'} \ldots i_{X_k'}\rho$. We claim that
\[\varphi^{\ast}=(-1)^kd \rho^{\ast}\] We prove this by induction on $k$.
Namely let $\varphi_l^{\ast}=i_{X_1'}\ldots i_{X_l'}\varphi$ and
$\rho_l^{\ast}= i_{X_1'} \ldots i_{X_l'} \rho$. We claim that
$\varphi_l^{\ast}=(-1)^ld\rho_l^{\ast}$. For $l=1$ we have that $\rho$ is
$X_1'$-invariant. Hence
\[0={\cal L}_{X_1'}\rho=d \rho_1^{\ast}+ \varphi_1^{\ast} \]
Now we use induction. The form $\rho_{l-1}^{\ast}$ is $X_l'$-invariant. Hence
\[ 0 = {\cal L}_{X_l'}\rho_{l-1}^{\ast}=d \rho_l^{\ast}+ 
(-1)^{l-1}\varphi_l^{\ast} \] and we are done by induction. 

Both $\varphi^{\ast}$ and $\rho^{\ast}$ are $T$-invariant. Let 
$\nu: Z \mapsto N_{red}$ and $\nu': Z'' \mapsto Q$ be the
quotient maps.
One easily sees that there is a unique $(n-k,0)$-form $\varphi'$ on $Q$ and
a unique $(n-k-1,0)$-form $\rho'$ on $Q$ s.t. \[
\nu^{\ast}(\varphi')=\varphi^{\ast} ~ , ~ \nu^{\ast}(\rho')= \rho^{\ast} ~ , ~
 \varphi'=(-1)^k d \rho' \]
We can define on $Q$ SLag submanifolds $L''$ by the conditions $\omega'|_{L''}=0 ~ , ~ Im \varphi'|_{L''}=0$. 

The vector field $Y$ is tangent to $Z''$ and $T$-invariant, hence it
projects to a vector field $Y'$ on $Q$. We had $i_Y \varphi=\rho$ on $K(N)$. Hence we also have $i_{Y'}\varphi'=\rho'$ on $Q$. We obviously have the following:
\begin{lem}
Let $L''$ be a SLag submanifold of $Q$, invariant under $Y'$. Then
$L'=\nu'^{-1}(L'')$ is a SLag submanifold of $K(N)$, invariant under $T^k$
and under the $Y$-flow.  
\end{lem}
The proof of the lemma is obvious. 
\subsection{Complexity one actions and periodic orbits}
In the previous section we have shown that one can reduce the problem of finding $T$ and $Y$-invariant SLag submanifolds of $K(N)$ to finding $Y'$-invariant
SLag submanifolds of $Q$. In this section we will assume that  $k=n-1$. 
Let $X \subset Z''$ be the set of elements in $Z''$ of unit length and $S=X/T
\subset Q$. We will show that there is a vector field $W$ on $S$
s.t. there is correspondence between $Y'$-invariant SLag submanifolds of
$Q$ and the trajectories of the $W$-flow on $S$. 

As we saw the tangent bundle of $K(N)$ decomposes as a direct sum $V \oplus H$ of the vertical and the horizontal distributions. Let $U$ be the image of the Lie algebra of $T$ under the differential of the action on $K(N)$. At points of $Z''$ $U$ is an $(n-1)$-dimensional vector space, and it is contained in the horizontal distribution $H$ (since on $\pi(Z'')$ the moment map $\mu$ vanishes). Also the Kahler form $\omega_u$ restricts to $0$ on $U$. Let $U^c$ be the complexification of $U$ in the tangent bundle to $K(N)$. Then $U^c$ can be viewed as a complex $(n-1)$-dimensional vector bundle over $Z''$. Let $H'$ be the orthogonal complement of $U^c$ in $H$. Then the tangent bundle of $Z''$ is a direct sum $V \oplus H' \oplus U$. Also the quotient of $V \oplus H'$ under the $T$-action can be identified with the tangent bundle to the symplectic reduction $Q=Z''/T$. Since $H'$ and $V$ are $T$-invariant the tangent bundle to $Q$ splits as a direct sum of 2 complex line bundles: $TQ=V \oplus H'$. Also $V$ and $H'$ are orthogonal both with respect to the symplectic form $\omega'$ and the Riemannian metric on $Q$.

There is a natural circle action on $X$ (given by the multiplication by complex numbers of absolute value 1 on $K(N)$). This action is $T$-invariant, hence it induces a circle action on $S=X/T$. Let $F$ be the vector field generating this action on $S$. Then $F =J(Y')$ (here $J$ is the complex structure on $Q$). Also both $Y'$ and $F$ are in the vertical distribution $V$ along $S$ and the tangent bundle $TS$ of $S$ splits as a direct sum $TS=H' \oplus span(F)$.

Let $\gamma$ be some path in $S$ and let $\gamma^Q$ be the orbit of $\gamma$ under the $Y'$-flow in $Q$. We wish to understand when $\gamma^Q$ is a SLag submanifold of $Q$. Let $W$ be a tangent vector to $\gamma$. Clearly for $\gamma^Q$ to be Lagrangian we need $W$ to live in the horizontal distribution $H'$. The form $\rho'=i_{Y'}\varphi'$ is a (non-zero) $(1,0)$-form on $H'$. Hence the form $Im \rho'$ has a 1-dimensional kernel in $H'$. Clearly for $\gamma^Q$ to be Special we need $W$ to belong to this kernel. We can normalize $W$ s.t. $Re\rho'(W)=1$. Those conditions give rise to a non-vanishing horizontal vector field $W$ on $S$. Let $\gamma$ be a trajectory of $W$ on $S$ and consider $\gamma^Q \subset Q$. The forms $\omega'$ and $\varphi'$ vanish on $\gamma^Q$ along $\gamma$. Also the $Y'$-low preserves the horizontal distribution and ${\cal L}_{Y'}\rho'=\rho'$. From this we easily deduce that $\gamma^Q$ is a $Y'$-invariant SLag submanifold of $Q$. From the above discussion we get the following:
\begin{lem}
Let $\gamma$ be a trajectory of $W$ on $S$. Then $L_{\gamma}= \pi(\nu'^{-1}(\gamma^Q))$ is an immersed minimal Lagrangian submanifold of $N$. If $\gamma$ is periodic then $L_{\gamma}$ is an immersed minimal Lagrangian torus.
\end{lem}

There is one general relation among trajectories of $W$, which will later be important: Consider the circle action on $K(N)$ as before. The $(n,0)$-form $\rho$ is equivariant with respect to this action, i.e. if $\lambda \in S^1$ then $\lambda^{\ast}(\rho)=\lambda \rho$. Since $\varphi=d\rho$ we get that $\varphi$ is also equivariant with respect to this action. Thus we also deduce that $\rho'$ and $\varphi'$ are equivariant with respect to the circle action on $Q$. Also this action preserves the horizontal distribution $H'$ on $S$. Consider an element $-1 \in S^1$. Then $-1^{\ast}(\rho')=-\rho'$. From this we deduce that the $-1$-action on $S$ reverses the vector field $W$, i.e. $-1_{\ast}(W)=-W$. Thus the $-1$-action sends $W$-trajectories to $W$-trajectories, but it reverses their directions. 

\section{Toric K-E manifolds}
In Section 2.3 we saw that if we have a $T^{n-1}$-action on $N$, then one can construct minimal Lagrangian submanifolds of $N$ from trajectories of the vector field $W$ on $S$. In order to do this we needed $T$ to act freely on $Z''$. In this section we will show a class of examples where this holds. We will also investigate periodic orbits of $W$ on $S$ (to construct immersed minimal Lagrangian tori).  

Let $N$ be toric, i.e. we have an effective structure-preserving $T^n$-action on $N$. For recent results on toric K-E manifolds we refer the reader to \cite{Tian} and \cite{Bat}.
We will use various $(n-1)$-dimensional sub-tori of $T$ to construct invariant minimal Lagrangian submanifolds. But first we will see that there is a unique minimal Lagrangian torus, invariant under the whole of $T$.

Suppose $L$ is a regular orbit of the $T$-action, which is a minimal submanifold. Then $L^K$ is a SLag submanifold of $K(N)$. The moment map $\mu'$ of $K(N)$ is constant on $L^K$. As we have seen in Section 2.2, we must have $\mu=0$ on $L$ i.e. $L \subset \mu^{-1}(0)$. By Atyiah's result \cite{At}, $\mu^{-1}(0)$ is connected, hence $L=\mu^{-1}(0)$. So if a regular orbit, which is a minimal submanifold, exists, it must coincide with $\mu^{-1}(0)$. Next we prove that such an orbit does exist.  
\begin{lem}
Let $(M^{2n},\omega)$ be a compact symplectic manifold and $g$ some Riemannian metric on $M$. Suppose that we have an effective Hamiltonian n-torus action on $M$, which preserves $g$. Then there is a regular orbit of the action, which is a minimal submanifold with respect to $g$. In fact this orbit maximizes volume among the orbits.
\end{lem}
{\bf Proof:} We have a moment map $\mu$ and smooth orbits are levels set of the moment map. For a regular orbit $L$ to be a minimal submanifold, it is obviously necessary to be a critical point for the volume functional on the orbits. We note that it is also sufficient. Indeed let $v$ be any element of the Lie algebra $\mathcal{T}$ of the torus $T^n$. Then $\mu(v)$ is $T^n$-invariant, and so is the gradient
$\nabla \mu(v)$. Also this gradient is orthogonal to the orbits. Consider now this gradient flow. It commutes with the $T^n$-action, hence it sends orbits to orbits. Since $L$ is critical for the volume functional on the orbits, we get from the first variation formula $\int_{L}h \cdot \nabla \mu(v)=0$. Here $h$ is a trace of the second fundamental form of $L$. But both $h$ and $\nabla \mu(v)$ are $T^n$-invariant, hence we are integrating a constant. So $h \cdot \nabla \mu(v)=0$ pointwise. Now $v$ was arbitrary, hence $h=0$.  

We want to find a regular orbit, which is maximum point for the volume functional on the orbits. First we need to prove that the volume functional is continuous on the space of orbits. Let $L'$ be a regular orbit for the torus action. Then the differential of the moment map is surjective along $L'$. From this one easily deduces that orbits of the action near $L'$ coincide with level sets of the moment map. So obviously the volume functional is continuous near $L'$. Next we prove that the volume functional is continuous near the singular orbits. This follows from the following easy Lemma:
\begin{lem}
Let $L$ be an orbit with a positive dimensional stabilizer $T' \subset T$ and $x \in L$. Then for any $\epsilon > 0$ there is a neighbourhood $U$ of $x$ s.t. any orbit passing through $U$ has volume $< \epsilon$.
\end{lem}
{\bf Proof:} Take a (unit) vector $e_1$ in the Lie Algebra of $T'$. Then the corresponding flow vector field $X_1$ vanishes along $L$. Extend $e_1$ to an o.n. basis $e_2,\ldots,e_n$ of ${\cal T}$. The flow vector fields $X_i$ will have uniformly bounded lengths. We choose a neighbourhood $U$ of $x$ in which $X_1$ has sufficiently small length and it is clear that volumes of the orbits through $U$ will be sufficiently small. Q.E.D.

So the volume functional is continuous on the space of orbits and we can find an orbit $L$, which maximizes volume among the orbits. Obviously $L$ must be a regular orbit (since singular orbits have zero volume). As we have seen, $L$ is a minimal submanifold of $N$ and we are done.           Q.E.D.

Let now $T'' \subset T^n$ be some $(n-1)$-torus in $T$ and let $\mu''$ be the canonical moment map for the $T''$-action on $N$ as in Section 2.2. Then $\mu''$ is just the restriction of $\mu$ to the dual Lie algebra of $T''$. In order to apply the constructions of Section 2.3 we want $T''$ to act freely on $Z''$. The following lemma guarantees the existence of such $T''$:
\begin{lem}
Let $N^{2n}$ be a K-E manifold with an effective $T^n$-action as above. Then there is an $(n-1)$- torus $T'' \subset T$ s.t.

i) The differential of the $T''$-action on $N$ is injective along $Z$ and
$T$ acts freely on $Z''$.

ii) There is an element $v$ in the Lie algebra of $T''$ s.t. the flow vector field $X_v$ doesn't have a constant length along $Z$.
\end{lem}
{\bf Remark:} Condition ii) in the lemma will be used to show that certain minimal Lagrangian tori we shall construct have Killing fields of non-constant length, hence they are not flat.

{\bf Proof:} Let $T'' \subset T^n$ be some $(n-1)$-torus. First we prove that if the differential of the $T''$-action on $N$ is injective along $Z$, then the $T$-action on $Z''$ is free. Suppose not. Then there is a point $l \in Z''$ and an element $1 \neq t \in T$ s.t $t \cdot l=l$. In that case $t$ also preserves the points on the $T''$-orbit through $l$. The tangent space $P$ to this orbit at $l$ is in the horizontal distribution at $l$ (since we are at the zero set of the moment map $\mu'$). Also $\omega_u|_P=0$. So the differential $dt$ of the $t$-action at $l$ acts trivially on the complexification $P^c$ of $P$. Also $dt$ acts trivially on the vertical distribution $V(l)$ at $l$. The vector space $P^c \oplus V(l)$ is a complex vector space of dimension $n$ and $dt$ acts trivially on it. Also $dt$ preserves the holomorphic volume form $\varphi$ at $l$. Hence $dt$ is trivial at $l$. Hence $t$ acts trivially on $K(N)$ and on $N$, but the $T$-action on $N$ was effective- a contradiction.

Next we wish to understand for which $(n-1)$-tori $T'' \subset T$ the differential of the $T''$-action is injective along the zero set $Z$ of the canonical moment map of $T''$. Let ${\cal T}^{\ast}$ be the dual Lie algebra of $T$ and let $\Lambda \subset {\cal T}^{\ast}$ be the weight lattice of $T$. Any element $0 \neq v \in \Lambda$ defines an $(n-1)$-torus $T_v \subset T$ s.t. $v$ vanishes on the Lie algebra of $T_v$. Let $\mu$ be the canonical moment map of $T$ and $\mu_v$ be the canonical moment map of $T_v$. Then $\mu_v$ is just the restriction of $\mu$ to the dual Lie algebra of $T$. It is therefore clear that $\mu_v$ vanishes at a point $n \in N$ iff $\mu(n)$ is proportional to $v$. Since $N$ is a toric variety, the moment polytope is convex and has no faces in the interior. Since $0$ is in the interior of the moment polytope, it is clear that $Z =\mu^{-1}[t_1v,t_2v]$ with $t_1 < 0 < t_2$. For any $t_1 < t < t_2$ the value $tv$ is in the interior of the moment polytope, while $t_1v$ and $t_2v$ are not. 

Suppose the line $span(v)$ doesn't intersect any of the $(n-2)$-faces of the moment polytope. This means that any point in $Z$ has either a trivial or a 1-dimensional stabilizer in $T$. We claim that the differential of the $T''$-action is injective along $Z$. Suppose not. Then there is a point $n \in Z$ and a vector $0 \neq w $ in the Lie algebra of $T''$ s.t. the flow vector field $X_w$ vanishes at $n$. Since $n \in Z$ the flow vector field $X_w'$ of $w$ on $K(N)$ is horizontal along $\pi^{-1}(n) \subset K(N)$, hence it vanishes along $\pi^{-1}(n)$. Let $g= exp(tw)$ for some $t \in \mathbb{R}$. Then the $g$-action on $\pi^{-1}(n)$ is trivial. But this means that the differential $dg$ of the $g$-action on the tangent space  $T_nN$ has Jacobian 1. Also $g$ acts trivially on the orbit $L'$ of the $T$-action through $n$. The tangent space $T_nL'$ of $L'$ at $n$ is $(n-1)$-dimensional and $\omega$ restricts to $0$ on it. Hence it's complexification $T_nL'^c$ is a complex $(n-1)$-dimensional space and $dg$ acts trivially on it. Also $dg$ has Jacobian 1. Hence $dg$ is trivial, hence $g$ acts trivially- a contradiction.

A generic line in the projective space $P{\cal T}^{\ast}$ doesn't intersect the $(n-2)$-faces of the moment polytope of $\mu$. Also the set of lines passing through points of $\Lambda$ is dense in $P{\cal T}^{\ast}$. So we can easily find $v \in \Lambda$ so that i) holds for $T_v$. In order to ensure that ii) holds, consider a point $b$ in the $(n-2)$-face of the moment polytope. The orbit $\mu^{-1}(b)$ has a stabilizer of dimension at least 2. Hence we can find a vector $0 \neq w \in {\cal T}$ s.t. $b(w)=0$ and the flow vector field $X_w$ vanishes along $\mu^{-1}(b)$. We can find a sequence of elements $v_k \in \Lambda$ s.t. the lines $(v_k)=span(v_k)$ do not intersect the $(n-2)$-faces of the moment polytope and $(v_k)$ converge to the line $(b)=span(b)$ in $P{\cal T}^{\ast}$. We can also find a sequence of vectors $w_k \in {\cal T}$
 s.t. $v_k(w_k)=0$ and $w_k \mapsto w$.  

Each $v_k$ defines an $(n-1)$-torus $T_k \subset T$. Let $\mu_k$ be the canonical moment map of $T_k$, and $Z_k$ be the zero set of $\mu_k$. We can find points $n_k$ on $Z_k$ s.t. $n_k$ converge to a point $n \in \mu^{-1}(b)$. Let $X_k$ be the flow vector field of $w_k$. Then the length of $X_k$ at points $n_k$ goes to $0$ as $k \mapsto \infty$. On the other hand the torus $L=\mu^{-1}(0)$ is contained in all of $Z_k$. Moreover the lengths of $X_k$ along $L$ are a-priori bounded from below. So we deduce that for $k$ large enough the torus $T'' = T_k$ satisfies the conditions i) and ii) of the lemma.  Q.E.D.

From now on we pick a sub-torus $T'' \subset T$ satisfying the conditions of Lemma 9. We can use the results of Section 2.3 to deduce that one can construct minimal Lagrangian submanifolds of $N$ from the trajectories of the vector field $W$ on $S$. From Lemma 6 we deduce that in order to obtain immersed minimal Lagrangian tori we need the orbits to be periodic. A first step in finding such orbits will be the following observation: The circle $R=T/T''$ acts freely on $Q$ and on $S$. Let $w \neq 0$ be some element in the Lie algebra of $R$. We have the flow vectors field $A_w$ for the $w$-action on $Q$ and the vector fields $A_w$ and $W$ commute. We also have a $(1,0)$-form $\rho'$ and a holomorphic $(2,0)$-form $\varphi'$ on $Q$ with $\varphi'=(-1)^{n-1}d\rho'$. The flow of $A_w$ preserves $\rho'$ and $\varphi'$.
A key point in finding periodic trajectories of $W$ is the fact that there is a function on $S$ constant along the trajectories:
\begin{lem}
Let $h=\rho'(A_w)$ and $f=Re(h)$. Then $f$ is constant along the trajectories of $W$.
\end{lem}
{\bf Proof:} We have \[0 = {\cal L}_{A_w}\rho'=d(i_{A_w}\rho') + i_{A_w}d\rho'=dh+(-1)^{n-1}i_{A_w}\varphi' \]
So $dh=(-1)^ni_{A_w}\varphi'$. So $dh(W)= (-1)^n\varphi'(A_w,W)$. The vector field $A_w$ is in the tangent bundle to $S$, hence we can decompose it into $A_w=A_w^H + \lambda F$. Here $A_w^H$ is the horizontal part of $A_w$ (i.e. the part in the distribution $H'$), $F$ is the generator of the $S^1$-action on $S$ and $\lambda \in \mathbb{R}$ (see Section 2.3). $W$ is horizontal and $H'$ is a 1-dimensional complex vector bundle. The form $\varphi'$ is a $(2,0)$-form on $Q$. Hence $\varphi'(A_w^H,W)=0$. Also $F=JY'$. Hence $\varphi' (F,W)= i\varphi'(Y',W)$. By the construction of $W$ we had that $\varphi'(Y',W)$ is real. From all this we deduce that $dh(W)$ is purely imaginary, hence $df(W)=0$, i.e. $f$ is constant along the trajectories of $W$. Q.E.D. 

From the previous lemma we deduce that the trajectories of $W$ live on level sets of the function $f$. We need to understand those level sets in more detail.

We had our symplectic reductions $N_{red}=Z/T''$ and $Q$ and we have a natural projection $\pi':Q \mapsto N_{red}$. Let $v$ be an element of the weight lattice $\Lambda$ of ${\cal T}^{\ast}$ defining the torus $T''$. As we have seen $Z$ is equal to $\mu^{-1}[t_1v,t_2v]$ for $t_1 < 0 < t_2$. $T''$ acts freely on $Z_0=\mu^{-1}(t_1v,t_2v)$ and we have $N_0=Z_0/T'' \subset N_{red}$, which is the smooth part of $N_{red}$. We also have 2 points $a_i = \mu^{-1}(t_iv)/T'' \in N_{red}$ and $N_{red}$ is a disjoint union of $a_1,a_2$ and $N_0$.  We have $S_0=\pi'^{-1}(N_0) \bigcap S$, and $S_0$ is a fiber bundle over $N_0$ with fibers being the orbits of the $S^1$-action on $S$ (see Section 2.3). This action is free on $S_0$. Also each $K_i =\pi'^{-1}(a_i)$ is an orbit of the $S^1$-action on $S$, but this action on each $K_i$ might have a finite stabilizer.

We have seen in Section 2.3 that the form $\rho'$ is equivariant with respect to the $S^1$-action on $S$. The flow vector field $A_w$ is invariant under $S^1$-action. Hence the function $h=\rho'(A_w)$ is $S^1$-equivariant.

On $Z_0$ we had an oriented Lagrangian distribution $D$, given by the image of ${\cal T}$ under the differential of the action on $N$. This distribution gives rise to a unit length section $\kappa$ of $K(N)$ over $Z_0$ as in Lemma 1. This section is $T$-invariant, hence it gives rise to an $R$-invariant section $\kappa'$ of $S_0$ over $N_0$. By definition $\kappa$ restricts to a positive real $n$-form on the distribution $D$. From this we deduce that $h=\rho'(A_w)$ is real and positive along $\kappa'$. 

We can normalize $w$ s.t. $v(w)=1$. We have a function $\tau=\mu(w)$ on $N_{red}$, and the image of this function is the interval $[t_1,t_2]$. For each $t_1 \leq t \leq t_2$ the level set $\tau^{-1}(t)$ is an orbit of the $R$-action on $N_{red}$. Let $L'=\tau^{-1}(0)$, $L_+=\kappa'(L') \subset S$ and $L_-= (-1) \cdot L_+$. Each $L_{\pm}$ is an orbit of the $R$-action on $S$. Also at points of $L_{\pm}$ the vector field $A_w$ is horizontal (since $\mu(w)=0$) and $\rho'(A_w)$ is real.
The vector field $W$ also satisfies those properties, hence $W$ is proportional to $A_w$ along $L_{\pm}$. So we see that $L_{\pm}$ are trajectories $W$ (of course the minimal Lagrangian torus of $N$ coming from these trajectories is the torus $L=\mu^{-1}(0)$). We have the following:
\begin{lem}
The differential $df$ of $f$ is non-vanishing on $S-(L_- \bigcup L_+)$.
\end{lem}
{\bf Proof:} We have seen in the proof of Lemma 10 that $dh=(-1)^n i_{A_w}\varphi'$. On $S \bigcap \pi'^{-1}(N_{red}-L')$ the vertical part of the vector field $A_w$ doesn't vanish. Hence the form $i_{A_w}\varphi'$ restricts as a non-vanishing $(1,0)$-form on the horizontal distribution $H'$. From this it is clear that $df|_{H'} \neq 0$.

On $S \bigcap \pi'^{-1}(L')-(L_- \bigcup L_+)$ $h$ is not real. Also $h$ is equivariant with respect to the $S^1$-action. Let $F$ be the vector field generating the $S^1$-action as before. Then the derivative of $f=Reh$ is non-zero in the direction of $F$.  Q.E.D.

$f$ attains a constant value $f_+$ along $L_+$ and a value $f_-=-f_+$ along $L_- = (-1) \cdot L_+$. Since $S$ is compact and connected, it is clear from Lemma 11 that $f_+$ is the absolute maximum of $f$, attained only at $L_+$, and $f_-$ is the absolute minimum of $f$, attained only at $L_-$. Also for any $s \in (f_-,f_+)$, the level set $\Sigma_s=f^{-1}(s)$ is smooth.

We will also need the fact that $f|_{K_i}=0$. To prove that we note that along $K_i$ $A_w$ is vertical. Indeed the action of $exp(tw)$ on $a_i$ is trivial for any $t \in \mathbb{R}$. Hence the action of $exp(tw)$ on $S$ preserves the fiber $K_i=\pi'^{-1}(a_i)$, so the vector field $A_w$ is tangent to $K_i$, i.e. vertical along $K_i$. But from this we deduce that $h=\rho'(A_w)=0$ at $K_i$, and so $f=0$ at $K_i$.

Let now $\Phi=f^{-1}(f_-,f_+)$. Take any point $m \in \Phi$ and consider the level set $\Sigma_s$ of $f$ passing through $m$. Let $\Sigma_s^0$ be the connected component of $\Sigma_s$ containing $m$. The vector field $W$ is tangent to $\Sigma_s^0$. We have a free $R$-action on $\Sigma_s^0$, and this action preserves the vector field $W$. Also $A_w$ is transversal to $W$ at all points of $\Sigma_s^0$. Indeed let $m' \in \Sigma_s^0$. If $m' \in S \bigcap \pi'^{-1}(N_{red}-L')$, then the vector field $A_w$ is not horizontal, so it can't be proportional to $W$. If $m' \in S \bigcap \pi'^{-1}(L')- (L_- \bigcup L_+)$, then $h=\rho'(A_w)$ is not real, while $\rho'(W)$ is real. So again $A_w$ and $W$ can't be proportional. Thus we get that the quotient of $\Sigma_s^0$ by the $R$-action is a circle and $W$ projects to a non-vanishing vector field on it. From this we deduce that the $W$-trajectory starting at $m$ will intersect the $R$-orbit of $m$. Suppose it intersects this orbit for the first time at a point $\xi(m)m$, $\xi(m) \in R$. This gives rise to a well-defined function $\xi:\Phi \mapsto R$. Clearly $\xi$ is continuous, $R$-invariant and constant along the trajectories of $W$. Also we have seen in Section 2.3 that the $-1$-action on $S$ sends $W$-trajectories to $W$-trajectories in the reverse direction. From this we easily deduce that \[ \xi(-1 \cdot m)=\xi(m)^{-1} \]
Obviously the trajectory through $m$ is periodic iff $\xi(m)$ is a root of unity in $R$. Let $R'$ be the set of roots of unity in $R$. Since $\xi$ is continuous, $\xi^{-1}(R')$ will be everywhere dense in $\Phi$ unless $\xi$ assumes a constant value not in $R'$ on some open subset of $\Phi$. The next lemma shows that it is impossible:
\begin{lem}
Suppose that $\xi$ is constant on some open set $U \subset \Phi$. Then $\xi$ is equal to a constant $g$ on the whole of $\Phi$ and $g^2=1$. 
\end{lem}
{\bf Proof:} Let $S_+=f^{-1}(0,f_+)$, $S_-=f^{-1}(f_-,0)$. Thus $S_-= -1 \cdot S_+$. Suppose that $\xi$ is constant on some open set $U \in \Phi$. Then $\xi$ is constant on some open ball $U'$ either in $S_+$ or in $S_-$. We can assume w.l.o.g. that $U' \subset S_+$. We note that $S_+$ is connected. In fact $S_+$ is given by \[S_+=(\kappa'(n)e^{i \theta}| n \in N_0 ~ , ~ -\pi/2 < \theta < \pi/2 )-(L_+ \bigcup L_-)\]
First we prove that $\xi$ is a constant $g$ on $S_+$. Let $A_w^H$ be the horizontal part of the vector field $A_w$. Since $S_+ \subset \pi^{-1}(N_0)$ we deduce that $A_w^H$ doesn't vanish on $S_+$. We also note that $A_w^H$ cannot be proportional to $JW$. Indeed suppose that $A_w^H=\lambda JW$ for some $\lambda \in \mathbb{R}$ at some point $m \in S_+$. But then \[h(m)=\rho'(A_w)=\rho'(A_w^H)=i\lambda\rho'(W) \] 
So $h(m)$ is purely imaginary, hence $f(m)=0$ - a contradiction. Since both $A_w^H$ and $W$ lie in $H'$, which is a complex 1-dimensional distribution, we deduce that we can find a function $b :S_+ \mapsto \mathbb{R}$ s.t. the vector fields $W'=A_w^H +bJA_w^H$ and $W$ are proportional. Hence the trajectories of $W'$ and $W$ coincide. We will use $W'$ instead of $W$ to prove that $\xi$ is constant on $S_+$. Suppose that for a point $m \in S_+$ it takes time $t(m)$ for the $W'$-flow to hit the $R$-orbit of $m$. We have the following:

\begin{lem}
$\xi(m)=exp(t(m)w)$
\end{lem}
{\bf Proof:} Let $\gamma$ be the trajectory of $W'$ through $m$ and  $\gamma'=\pi'(\gamma)$ be the corresponding path in $N_0$. We have an $R$-action on $N_0$, and the corresponding flow vector field $B_w$ for the $w$-flow on $N_0$. We obviously have $\pi'_{\ast}(A_w^H)=B_w$. Hence the tangent field to $\gamma'$ is $B_w +bJB_w$.

The $R$-action on $N_0$ is Hamiltonian with the moment map $\tau=\mu(w)$. Also the $B_w$-flow on $N_0$ commutes with the complex structure $J$ on $N_0$. Hence the vector fields $B_w$ and $JB_w$ commute. Let $P_1=\gamma'(0)$ and $P_2=\gamma'(t(m))$. Then $P_2= \xi(m) P_1$. We need to prove that $P_2 =exp(t(m)w)P_1$ and since the $R$-action on $N_0$ is free we would be done.

Let $exp (xJB_w)$ be the time $x$ flow of $B_w$. Note that the $JB_w$-flow is not complete. In fact we have \[JB_w(\tau)=\omega_{red}(JB_w,B_w)=|B_w|^2>0 \]
So $\tau$ increases on the $JB_w$-trajectories. Let $c(r)= \int_{(0,r)}b(t)dt$ for $0 \leq r \leq t(m)$. Consider a path $\gamma''(r)= exp(c(r)(JB_w))exp(rw)(P_1)$ (note that we flow $P_1$ with respect to $rB_w$ first). Then $\gamma''(r)$ is defined for small values of $r$. Also the tangent vector to $\gamma''$ is $B_w+b(r)JB_w$. So $\gamma''$ coincides with $\gamma'$ whenever it is defined. 

Suppose on $\gamma'$ $\tau$ ranges between $s_1$ and $s_2$. Then $t_1 < s_1$ and $s_2 <t_2$. Pick any $r$ for which $\gamma''(r)$ is defined and consider the path $exp(tJB_w)exp(rw)(P_1)$ for $t$ ranging between $0$ and $c(r)$. The function $\tau$ is increasing along the path, and on the endpoints it's values are between $s_1$ and $s_2$. Hence this path lives in the compact set $A=\tau^{-1}[s_1,s_2]$ in $N_0$. From this one can easily deduce that $\gamma''(r)$ is well defined for all $0 \leq r \leq s$ and coincides with $\gamma'(r)$. In particular $P_2=exp(c(t(m))JB_w)exp(t(m)w)(P_1)$. Now \[ \tau(P_2)=\tau(P_1)=\tau(exp(t(m)w)P_1) \] and $\tau$ increases on the trajectories of $JB_w$. So we get that $c(t(m))=0$, i.e. $P_2=exp(t(m)w)P_1$.  Q.E.D.

Now we can prove that $\xi$ is constant on $S_+$. Since $\xi$ is constant on $U'$ we get that $t(m)$ is a constant $t$ on $U'$. Let $\phi_t$ be the time $t$ flow of $W'$ on $S_+$. Consider the map $\chi=exp(-tw) \cdot \phi_t: S_+ \mapsto S_+$. $S_+$ is a connected real analytic manifold and $\chi$ is a real analytic map. Also $\chi$ is the identity map on $U'$. So we deduce that $\chi$ is the identity map. So $\phi_t$ is the multiplication by $g=exp(tw)$ on $S_+$. From this we easily deduce that $\xi=g$ on $S_+$. 

So $\xi$ assumes a constant value $g$ on $S_+$, and hence it assumes a constant value $g^{-1}$ on $S_-= -1 \cdot S_+$. Let $\Delta=f^{-1}(0)$. Then $\Delta$ is the common boundary of $S_+$ and $S_-$ in $\Phi$. Since $\xi$ is continuous, we must have $g=g^{-1}$, i.e. $g^2=1$.    Q.E.D.

From this we get an immediate corollary
\begin{cor}
The set $\xi^{-1}(R')$ is everywhere dense in $\Phi$.
\end{cor}
We are now ready to state and prove our main result:
\begin{thm}
Let $N^{2n}$ be a K-E manifold with positive scalar curvature with an
effective $T^n$-action. Then precisely one regular orbit of the action is a
minimal Lagrangian submanifold of $N$. Moreover there is an $(n-1)$-torus
$T^{n-1} \subset T^n$ and a sequence of non-flat $T^{n-1}$-invariant immersed
minimal Lagrangian tori $L_k \subset N$ s.t. $L_k$ locally converge to $L$ (in particular the supremum of sectional curvatures of
$L_k$ and the distance between $L$ and $L_k$ goes to $0$ as $k \mapsto \infty$). 
\end{thm}
{\bf Proof:} Choose a torus $T''=T^{n-1}$ which satisfies the conditions of Lemma 9. By Corollary 1 we can choose a sequence of points $m_k \in \xi^{-1}(R')$ s.t. $m_k$ converge to a point $m$ in $L_+$. The $W$-trajectories $\gamma_k$ through $m_k$ are periodic and live on level sets $\Sigma_k$ of $f$ with $\Sigma_k$ converging to $L_+$. From this we easily see that $\gamma_k$ locally converge to the trajectory $L_+$. One easily deduces that the immersed minimal Lagrangian tori $L_k$ which $\gamma_k$ define as in Lemma 6 locally converge to the minimal, $T$-invariant orbit $L$. 

Finally we prove that $L_k$ are not flat. From Lemma 9 we get a vector $v$ in the Lie algebra of $T''$ s.t. the flow vector field $X_v$ of $v$ doesn't have a constant length on $Z$. Now the vector field $X_v$ along $L_k$ is a Killing vector field of $L_k$. So to prove that $L_k$ is not flat it is enough to prove that $|X_v|^2$ is non-constant on $L_k$.

The function $|X_v|^2$ is $T$-invariant on $Z$. Thus it can be viewed as an $R$-invariant function on $N_{red}$, i.e. it can be viewed as a function of $\tau=\mu(w)$ on $N_{red}$. Also $|X_v|^2$ is a real analytic function on $N_0=\tau^{-1}(t_1,t_2)$. Since $|X_v|^2$ is non-constant, it is nowhere a locally constant function of $\tau$. Since $\gamma_k$ are different from $L_{\pm}$, we easily deduce that $\tau(\pi'(\gamma_k))$ are non-trivial intervals in $(t_1,t_2)$. Hence $|X_v|^2$ is non-constant on $L_k$ and we are done.  Q.E.D.

Massachusetts Institute of Technology

E-Mail : egold@math.mit.edu

\end{document}